\documentclass[11pt]{article}
\usepackage[latin1]{inputenc}
\usepackage{a4}
\usepackage{amssymb}
\usepackage{amsmath}
\begin{document}
\setlength{\parindent}{1.2em}
\def\COMMENT#1{}
\def\TASK#1{}
\let\TASK=\footnote
\def\noproof{{\unskip\nobreak\hfill\penalty50\hskip2em\hbox{}\nobreak\hfill%
       $\square$\parfillskip=0pt\finalhyphendemerits=0\par}\goodbreak}
\def\endproof{\noproof\bigskip}
\newdimen\margin   
\def\textno#1&#2\par{%
   \margin=\hsize
   \advance\margin by -4\parindent
          \setbox1=\hbox{\sl#1}%
   \ifdim\wd1 < \margin
      $$\box1\eqno#2$$%
   \else
      \bigbreak
      \hbox to \hsize{\indent$\vcenter{\advance\hsize by -3\parindent
      \sl\noindent#1}\hfil#2$}%
      \bigbreak
   \fi}
\def\proof{\removelastskip\penalty55\medskip\noindent{\bf Proof. }}
\def\enddiscard{}
\long\def\discard#1\enddiscard{}
\newtheorem{firstthm}{Proposition}
\newtheorem{thm}[firstthm]{Theorem}
\newtheorem{prop}[firstthm]{Proposition}
\newtheorem{lemma}[firstthm]{Lemma}
\newtheorem{cor}[firstthm]{Corollary}
\newtheorem{problem}[firstthm]{Problem}
\newtheorem{defin}[firstthm]{Definition}
\newtheorem{conj}[firstthm]{Conjecture}
\def\eps{{\varepsilon}}
\def\N{\mathbb{N}}
\def\R{\mathbb{R}}
\def\cH{\mathcal{H}}
\newcommand{\ex}{\mathbb{E}}
\newcommand{\pr}{\mathbb{P}}

\title{Maximizing several cuts simultaneously}
\author{Daniela K\"uhn \and Deryk Osthus}
\date{}
\maketitle
\vspace{-.8cm}
\begin{abstract} \noindent
Consider two graphs $G_1$ and $G_2$ on the same vertex set $V$ and suppose that
$G_i$ has $m_i$ edges. Then there is a bipartition of $V$ into two classes 
$A$ and $B$ so that for both $i=1,2$ we have
$e_{G_i}(A,B) \ge m_i/2-\sqrt{m_i}$. This answers a question of Bollob\'as and Scott.
We also prove results about partitions into more than
two vertex classes.
 
\end{abstract}

\section{Introduction}
Given a graph $G$ with $m$ edges, the Max-Cut problem is to determine
(the size of) the maximum cut in 
$G$. For complete graphs, the largest cut has size $m/2+o(m)$.
On the other hand, it is well known that a cut of
size at least $m/2$ in a graph $G$ can be found using the natural greedy 
algorithm. Improving this, Edwards~\cite{Edwards73, Edwards75} showed that every graph
with $m$ edges has a cut of size
$$m/2+\sqrt{\frac{m}{8}+\frac{1}{64}}-\frac{1}{8},
$$
which is
best possible. The Max-Cut problem is equivalent to finding a bipartition
$V_1,V_2$ of the vertex set of $G$ which minimizes $e_G(V_1)+e_G(V_2)$, where
$e_G(V_i)$ denotes the number of edges in the subgraph of $G$ induced by~$V_i$.
The related problem when one is looking for a partition into $k$ classes $V_1,\dots,V_k$
which minimizes all $e_G(V_i)$ simultaneously, i.e. which minimizes
$\max\{e_G(V_1),\dots,e_G(V_k)\}$, was studied by Bollob\'as and Scott~\cite{BS93,BS99,BS_JGT}
as well as Porter~\cite{Porter92, Porter94, Porter99}, see also~\cite{BS02}
for a survey. 

Here, we suppose that we are given several graphs on the same vertex set
and we want to find a bipartition which maximizes the sizes of the cuts
for all these graphs simultaneously. This problem was posed by Bollob\'as
and Scott~\cite{BS_JGT}. More precisely, they asked the following
question: What is the largest integer $f(m)$ such that whenever $G_1$ and
$G_2$ are two graphs with $m$ edges on the same vertex set $V$, there exists
a bipartition of $V$ in which for both $i=1,2$ at least $f(m)$ edges of
$G_i$ go across (i.e.~their endvertices lie in different partition classes).
They suggested that perhaps even $f(m)=(1-o(1))m/2$, i.e.~that we can
almost do as well as in the case where we
only have a single graph.  
Theorem~\ref{main} shows that this is indeed the case.

Given a graph $G$ and disjoint subsets $A,B$ of its vertex set, let 
$e_G(A,B)$ denote the number of edges between $A$ and~$B$.  

\begin{thm} \label{main}
Consider graphs $G_1,\dots,G_\ell$ on the same vertex set $V$ and suppose that
$G_i$ has $m_i$ edges. Then there is a bipartition of $V$ into two classes 
$A$ and $B$ so that for all $i=1,\dots,\ell$ we have
$$e_{G_i}(A,B) \ge \frac{m_i}{2}-\sqrt{\ell m_i/2}.$$
\end{thm}

Rautenbach and Szigeti~\cite{RS} observed that even for $\ell=2$ we cannot
guarantee that $e_{G_i}(A,B) \ge m_i/2$ for all~$i$. Indeed, let $G_1$
and $G_2$ be two edge-disjoint cycles of length~5 on the same vertex set.
(So $G_1\cup G_2=K_5$.) They also proved that $f(m)\ge m/2-\Delta^3$
if $\Delta(G_i)\le \Delta$ for $i=1,2$. (This answers the problem of
Bollob\'as and Scott if $(\Delta(G_i))^3=o(m)$ for $i=1,2$.)	

The following result for partitions of graphs into more than 
two parts shows that simultaneously for all graphs we can ensure that the number
of crossing edges is almost as large as one would expect in a random partition
(and almost the value one can ensure if one partitions only a single graph).

\begin{thm} \label{kpartite}
Let $k\ge 2$.
Consider graphs $G_1,\dots,G_\ell$ on the same vertex set $V$ and suppose that
$G_i$ has $m_i$ edges. Then there is a partition of $V$ into $k$ classes  $V_1,\dots,V_k$
so that for all $i=1,\dots,\ell$ the number of edges spanned by the
$k$-partite subgraph of $G_i$ induced by $V_1,\dots,V_k$ is at least 
$$\frac{(k-1)m_i}{k}-\sqrt{2\ell m_i}.$$ 
\end{thm}
 
In fact, if $\Delta(G_i)=o(m_i)$ for each $i$, then we can strengthen the conclusion:
The next theorem shows that there is a partition of $V$ into $k$ classes where each
of the $\binom{k}{2}$ bipartite graphs spanned by two of the partition classes contains almost
$2m_i/k^2$ edges for all $i=1,\dots,\ell$ simultaneously. Again, this is
about the number of edges which one would expect in a random partition.

\begin{thm} \label{kpartite2}
Let $k\ge 2$ and $0<\eps\le 1/(9\ell^2 k^4)$.
Consider graphs $G_1,\dots,G_\ell$ on the same vertex set $V$. Suppose that
$G_i$ has $m_i$ edges and that $\Delta(G_i) \le \eps m_i$ for all $i=1,\dots,\ell$.
Then there is a partition of $V$ into $k$ classes $V_1,\dots,V_k$
so that for all $i=1,\dots,\ell$ and for all $s,t$ with $1\le s < t\le k$
we have
$$e_{G_i}(V_s,V_t)\ge \frac{2m_i}{k^2}-\eps^{1/4} m_i$$
and $$e_{G_i}(V_s)\ge \frac{m_i}{k^2}-\eps^{1/4}m_i.$$ 
\end{thm}

Note that even for $\ell=1$ the condition that $\Delta(G_i)\le \eps m_i$
cannot be omitted completely.
For example, the result is obviously false if $G$ is a star.
On the other hand, a result of Bollob\'as and
Scott~\cite[Thm.~3.2]{BS_JGT} implies that in the case when the maximum degree of
each $G_i$ is bounded by a constant~$\Delta$, the bound on $e_{G_i}(V_s,V_t)$ in
Theorem~\ref{kpartite2} can be improved to $2m_i/k^2-C$ where $C=C(\ell,\Delta)$
(and similarly for $e_{G_i}(V_s)$). Note that this implies that if $G$ has bounded maximum 
degree, then one can achieve a bounded error term in Theorems~\ref{main} and~\ref{kpartite} as
well. 

The proofs of Theorems~\ref{main}--\ref{kpartite2} can be derandomized to
yield polynomial time algorithms which 
find the desired partitions (see Section~\ref{sec:alg}).  

\section{An open problem}
Consider an $r$-uniform hypergraph $\cH$ with $m$ hyperedges.
It is easy to see that there is a partition $V_1,\dots,V_r$ of the vertex
set of $\cH$ such that at least $r!m/r^r$ hyperedges of $\cH$ meet every~$V_i$
(in other words, each $r$-uniform hypergraph contains an $r$-partite subhypergraph with
at least $r!m/r^r$ hyperedges). To verify this, consider the expected number of hyperedges
which meet every~$V_i$ in a random partition of the vertices.
We believe that one does not loose much if one considers several
hypergraphs simultaneously:

\begin{conj}\label{hypergraphconj}
Suppose that $\cH_1,\dots,\cH_\ell$ are $r$-uniform hypergraphs on the same
vertex set~$V$ such that $\cH_i$ has $m_i$ hyperedges. Then there exists a partition
of $V$ into $r$ classes $V_1,\dots,V_r$ such that for all $i=1,\dots,\ell$ at
least $r!m_i/r^r-o(m_i)$ hyperedges of $\cH_i$ meet each of the classes $V_1,\dots,V_r$.
\end{conj}

Given an $r$-uniform hypergraph $\cH$ and distinct vertices $x,y\in \cH$, denote
by $N_{\cH}(x,y)$ the number of hyperedges which contain both $x$ and $y$.
Let $\Delta_2(\cH)$ denote the maximum of $|N_{\cH}(x,y)|$ over all pairs
$x\neq y$. 
One can adapt our proof of Theorem~\ref{kpartite2} to show that
Conjecture~\ref{hypergraphconj} holds in the case when $\Delta_2(\cH_i)=o(m_i)$ for each~$i$.
We omit the details.

\section{Proofs}

The proofs all proceed by considering a random partition and analyzing this using
the second moment method.

\begin{lemma}\label{randompart}
Let $c \in \R$ with $c > 1/2$. Suppose that $G$ is a graph with $m$ edges
whose vertex set is $V$.
Consider a random bipartition of $V$ into two classes $A$ and $B$ which is obtained
by including each $v \in V$ into $A$ with probability $1/2$ independently
of all other vertices in $V$. Then with probability at least $1-1/(2c)$ we have
$$e_{G}(A,B) \ge \frac{m}{2}-\sqrt{cm/2}.$$ 
\end{lemma}

If we apply the above result with $c=\ell$ (say) to the graphs in Theorem~\ref{main},
the failure probability for each of them is less than $1/(2\ell)$. Summing up
all these failure probabilities immediately implies Theorem~\ref{main}.

{\removelastskip\penalty55\medskip\noindent{\bf Proof of Lemma~\ref{randompart}. }}
For every edge $e$ of the graph $G$, define an indicator variable 
$X_e$ as follows: if one endvertex of $e$ is in $A$ and the other one is in $B$, then
let $X_e:=1$, otherwise let $X_e:=0$. Clearly, $\pr[X_e=1]=1/2$.
Also, for $e,e' \in E(G)$ with $e \neq e'$, we have
$$
\ex [X_e \cdot X_{e'}]= \pr[X_e=1,\, X_{e'}=1]= \frac{1}{2} \pr[X_e=1\mid X_{e'}=1]=\frac{1}{4}.
$$
Note that the final equality holds regardless of whether $e$ and $e'$ have an endvertex in 
common or not.
Now let $X:=\sum_{e \in E(G)} X_e$.
Thus $X$ counts the number of edges between $A$ and $B$ and $\ex X=m/2$.
Let $\sum_{e,e' \in E(G) \atop e \neq e'}$ denote the sum over all
ordered pairs $e,e'$ of distinct edges in~$G$. Then,
using the fact that $\ex[X_e^2]=\ex[X_e]$, we have
\begin{align*}
\ex[X^2] & = \sum_{e \in E(G)} \ex[X_e]+\sum_{e,e' \in E(G) \atop e \neq e'} \ex[X_e \cdot X_{e'}]\\
         & = \ex[X]+ \sum_{e,e' \in E(G) \atop e \neq e'} \frac{1}{4} = \frac{m}{2}+ \frac{m(m-1)}{4}
         = \frac{m(m+1)}{4}.
\end{align*}
This in turn implies that the variance of $X$ satisfies
${\rm Var} X=\ex[X^2]-(\ex X)^2=m/4$.
The result now follows from a straightforward application of Chebyshev's inequality:
$$
\pr[X \le m/2-\sqrt{c m/2}] \le \pr[ |X -\ex X| \ge \sqrt{c m/2}]
\le \frac{{2\rm Var}X}{c m}=\frac{1}{2c}.
$$
\endproof

\removelastskip\penalty55\medskip\noindent{\bf Proof of Theorem~\ref{kpartite}.}
As in Lemma~\ref{randompart}, we first consider a single graph $G$
with $m$ edges and vertex set~$V$. Consider a random partition of $V$ into
$k$ disjoint sets $V_j$ which is obtained by including each
$v\in V$ into $V_j$ with probability $1/k$ independently
of all other vertices. Let $X_e:=0$ if the edge $e$ has both its
endpoints in some $V_j$ and let $X_e:=1$ otherwise. So $\pr [X_e=1]=(k-1)/k$.
Also, it is easy to check that $\ex [X_e \cdot X_{e'}]=(k-1)^2/k^2$.
Again, this holds regardless of whether $e$ and $e'$ have an endvertex in 
common or not. Let $X$ denote the number of edges whose endvertices lie in different
vertex classes. Thus $\ex X=\frac{k-1}{k}m$ and
\begin{align*}
\ex[X^2] & = \sum_{e \in E(G)} \ex[X_e]+
\sum_{e,e' \in E(G) \atop e \neq e'} \ex[X_e \cdot X_{e'}]\\
& =\frac{k-1}{k}m+ m(m-1) \frac{(k-1)^2}{k^2}\le m+(\ex[X])^2.
\end{align*}
Therefore ${\rm Var} X\le m$ and so Chebyshev's inequality implies that
$$
\pr[X \le \frac{k-1}{k}m-\sqrt{2\ell m}] \le \pr[ |X -\ex X| \ge \sqrt{2\ell m}]
\le \frac{{\rm Var}X}{2\ell m}\le \frac{1}{2\ell}.
$$
Theorem~\ref{kpartite} now follows by summing up this bound on the failure
probability for each of the graphs $G_i$.
\endproof 
 
\removelastskip\penalty55\medskip\noindent{\bf Proof of Theorem~\ref{kpartite2}. }
Let $\eps$ be as in the statement of the theorem. 
As in the previous proof, we first consider a single graph $G$, this time with $m$ edges
and maximum degree $\Delta\le \eps m$.
Consider a random partition of $V:=V(G)$ into $k$ disjoint sets $V_j$ which is obtained
by including each vertex $v\in V$ into $V_j$ with probability $1/k$
independently of all other vertices. Fix some $s$ and $t$ with $1 \le s <t \le k$.
This time let $X_e:=1$ if one endvertex of $e$ is contained
in $X_s$ and the other in $X_t$. Put $X_e:=0$ otherwise. So $\pr[X_e=1]=2/k^2=:\alpha$.
Now the value of $\ex [X_e \cdot X_{e'}]$ depends on whether $e$ and $e'$ have an 
endvertex in common or not:
If they do have an endvertex in common, we will use the trivial bound
$\ex [X_e \cdot X_{e'}] \le 1 < 1+\alpha^2$. Note that the number of ordered
pairs $e, e'$ of distinct edges for which this can happen is trivially at most $2\Delta m$.
If $e$ and $e'$ have no vertex in common, then it is easy to see that
$$
\ex [X_e \cdot X_{e'}]=\pr[X_e=1]\pr[X_{e'}=1]=\alpha^2.
$$
Let $X:=\sum_{e \in E(G)} X_e$. Thus $\ex[X]=2m/k^2=\alpha m$. 
Moreover
\begin{align*}
\ex[X^2] & = \sum_{e \in E(G)} \ex[X_e]+\sum_{e,e' \in E(G) \atop e \neq e'} \ex[X_e \cdot X_{e'}]\\
         & < \ex[X]+ 2\Delta m+ \sum_{e,e' \in E(G) \atop e \neq e'} \alpha^2   \\
         & \le \alpha m+2\Delta m+ \alpha^2 m^2 \le 3\Delta m+(\ex[X])^2.
\end{align*}
Thus ${\rm Var}X \le 3\Delta m \le 3 \eps m^2$. So we can conclude that
$$
\pr[X \le \alpha m-\eps^{1/4} m] \le \pr[ |X -\ex X| \ge \eps^{1/4} m] 
\le \frac{{\rm Var}X}{\sqrt{\eps} m^2} \le 3\sqrt{\eps}\le \frac{1}{\ell k^2}.
$$
In exactly the same way one can show that $\pr[e_G(V_s)\le m/k^2-\eps^{1/4}m]\le 1/(\ell k^2)$.
(This time $\alpha:=1/k^2$.) 
Now sum up these failure probabilities for all the $\binom{k}{2}$ pairs $s,t$ and all
the $k$ values of $s$ to see that
the probability that a random partition does not have the required properties for $G$ is
at most $3/(4\ell)$. Again, Theorem~\ref{kpartite2} follows from summing up
this probability for all $G_i$.
\endproof 

We remark that at the expense of increasing the error terms the
partition classes in Theorems~\ref{main}--\ref{kpartite2} can be chosen to have
almost equal sizes. Indeed, Chernoff's inequality implies that in a
random partition of the vertex set as considered in the proofs with high
probability the vertex classes have almost equal sizes.

\section{Algorithmic aspects}\label{sec:alg}
Papadimitriou and Yannakakis~\cite{PY91} showed that the Max-Cut problem is APX-complete.
On the other hand,
as mentioned in the introduction, the obvious greedy algorithm always guarantees a 
cut whose size is at least $m/2$.
Moreover, the proofs described in the previous section can be derandomized 
to yield polynomial algorithms which construct partitions
satisfying the bounds in Theorems~\ref{main}--\ref{kpartite2}.
As the derandomization argument is similar for all three results, we only
only describe it for Theorem~\ref{main}.
More background information on derandomization can be found for instance 
in the books~\cite{ASp,MRbook} 
and in Fundia~\cite{Fundia} (in particular, the framework 
described in the latter applies to our situation). For simplicity, we consider
Theorem~\ref{main} only for~$\ell=2$, i.e.~in the case of two graphs.

So let $G_1$ and $G_2$ be two graphs whose vertex set is $V$
with $e(G_i)=m_i$. Consider a 
random partition of $V$ into sets $A$ and $B$ 
as described in the proof of Theorem~\ref{main} (cf.~Lemma~\ref{randompart}).
For $i=1,2$ define random variables $X_i:=e_{G_i}(A,B)$ and put
$\mu_i:=m_i/2=\ex[X_i]$. Set
$$Z_i:=\frac{\mu_i^2-2\mu_i X_i  +X_i^2}{m_i}
$$ for $i=1,2$ and
$Z:=Z_1+Z_2$. The proof of Theorem~\ref{main}
shows that for each $i$
$$\pr [  X_i < \mu_i -\sqrt{m_i}]\le \frac{{\rm Var}X_i}{m_i}<1/2.
$$ 
But $\ex[Z_i]={\rm Var}X_i/m_i$ and so
$\ex[Z]=\ex [Z_1]+\ex[Z_2]<1$.
Let $v_1,\dots,v_n$ be an enumeration of the vertices in~$V$.
Let $A_i$ denote the event that the vertex $v_i$
is contained in $A$. Then
$$
1>\ex[Z]=(\ex[Z \mid A_1]+\ex[Z \mid A_1^c] )/2 \ge \min 
\{\ex[Z \mid A_1],\ex[Z \mid A_1^c] \}.
$$
Thus at least one of $\ex[Z\mid A_1]$, $\ex[Z \mid A_1^c]$ has
to be less than~1. Let $C_1\in\{A_1,A_1^c\}$ be such that
$\ex[Z\mid C_1]<1$. Note that both $\ex[Z \mid A_1]$ and $\ex[Z \mid A_1^c]$
can be computed in polynomial time and so also $C_1$ can be determined in
polynomial time. Now
$$
1>\ex[Z\mid C_1]=(\ex[Z \mid C_1\cap A_2]+\ex[Z \mid C_1\cap A_2^c] )/2.
$$
So similarly as before there exists $C_2\in\{A_2,A_2^c\}$ such that
$\ex[Z\mid C_1\cap C_2]<1$ and $C_2$ can be determined in polynomial time.
We continue in this fashion until we have obtained events $C_k\in\{A_k,A_k^c\}$
for all $k=1,\dots,n$ such that
$$
\ex[Z\mid C_1\cap \dots \cap C_n]<1.
$$  
The proof of Chebyshev's inequality shows that
for each $i=1,2$ and for any event $U$ which has positive probability, we have
\begin{equation*} 
\pr [  X_i < \mu_i -\sqrt{m_i} \mid U] \le
\frac{\mu_i^2-2\mu_i\ex [X_i \mid U] +\ex[X_i^2 \mid U]}{m_i}=\ex[Z_i\mid U]
\end{equation*}
(the above also follows from Corollary~4 in~\cite{Fundia}).
Taking $U:=C_1\cap \dots \cap C_n$ this implies that
\begin{align}\label{eqfinal}
\sum_{i=1,2}\pr [  X_i< \mu_i-\sqrt{m_i} \mid U]\le \sum_{i=1,2}\ex[Z_i\mid U]
=\ex[Z\mid U]<1.
\end{align}
But $U:=C_1\cap \dots \cap C_n$ means that for each vertex $v_k\in V$ we
have decided whether $v_k\in A$ or $v_k\in B$. So the left hand side of~(\ref{eqfinal})
is either $0$ or~$1$, i.e.~it has to be~0. This means that the unique partition corresponding
to $C_1\cap \dots \cap C_n$ is as desired in Theorem~\ref{main}. Since each
$C_k$ can be determined in polynomial time this gives us a polynomial algorithm.

\section*{Acknowledgement}
We are grateful to Dieter Rautenbach for telling us about the problem.

{\footnotesize
\bigskip\obeylines\parindent=0pt
Daniela K\"uhn \& Deryk Osthus
School of Mathematics
Birmingham University
Edgbaston
Birmingham B15 2TT
UK
{\it E-mail addresses}: {\tt \{kuehn,osthus\}@maths.bham.ac.uk}
}

\end{document}